\documentclass[journal,twoside,web]{ieeecolor}
\usepackage{generic}
\usepackage{cite}
\usepackage{amsmath,amssymb,amsfonts}
\usepackage{algorithmic}
\usepackage{graphicx}
\usepackage{hyperref}
\hypersetup{hidelinks=true}
\usepackage{textcomp}

\usepackage{cite}
\usepackage{amsmath,amssymb,amsfonts,bm}
\usepackage{graphicx}
\usepackage{booktabs}
\usepackage{makecell}
\usepackage{xcolor}
\usepackage{url}
\usepackage[boxruled]{algorithm2e}
\usepackage[draft]{todonotes}

\usepackage{hyperref}

\newcounter{remarkcount}
\newcounter{assumpcount}
\newcounter{propcount}
\newcounter{lemmacount}
\newcommand{\lemma}[1]{
    \stepcounter{lemmacount}
    \textbf{Lemma} \arabic{lemmacount}: #1

}

\newcommand{\proposition}[1]{
    \stepcounter{propcount}
    {\textbf{Proposition} \arabic{propcount}}: #1
}

\newcommand{\remarknew}[1]{%
    \stepcounter{remarkcount}
    \textbf{Remark \arabic{remarkcount}:} #1
}

\newcommand{\assumption}[1]{
    \stepcounter{assumpcount}
    \textbf{Assumption \arabic{assumpcount}:} #1
}

\hypersetup{hidelinks} 

\newcommand{\revision}[1]{{\color{black} #1}}
\newcommand{\rmt}{\mathrm{T}}

\def\BibTeX{{\rm B\kern-.05em{\sc i\kern-.025em b}\kern-.08em
    T\kern-.1667em\lower.7ex\hbox{E}\kern-.125emX}}
\begin{document}
\title{Distributionally Robust System Level Synthesis With Output Feedback Affine Control Policy
}
\author{{Yun Li, Jicheng Shi, Colin N. Jones, Neil Yorke-Smith, Tamas Keviczky}
\thanks{Yun Li and Tamas Keviczky are with the Delft Center for Systems and Control, TU Delft, Delft, 2600 AA, the Netherlands (e-mail: y.li-39@tudelft.nl; T.Keviczky@tudelft.nl). }
\thanks{Jicheng Shi and Colin N. Jones are with the Automatic Control Laboratory, EPFL, 1015, Lausanne, Switzerland (e-mail: jicheng.shi@epfl.ch; colin.jones@epfl.ch).}
\thanks{Neil Yorke-Smith is with the STAR Lab, TU Delft, Delft, 2600 AA, the Netherlands (e-mail: n.yorke-smith@tudelft.nl).}
}

\maketitle

\begin{abstract}
This paper studies the finite-horizon robust optimal control of constrained linear systems subject to model mismatch and additive stochastic disturbances. Utilizing the system level synthesis (SLS) parameterization, we propose a novel SLS design using an output-feedback affine control policy and extend it to a distributionally robust setting to improve system resilience by minimizing the cost function while ensuring constraint satisfaction against the worst-case uncertainty distribution. The scopes of model mismatch and stochastic disturbances are quantified using the 1-norm and a Wasserstein metric-based ambiguity set, respectively. For the closed-loop dynamics, we analyze the distributional shift between the predicted output-input response -- computed using nominal parameters and empirical disturbance samples -- and the actual closed-loop distribution, highlighting its dependence on model mismatch and SLS parameterization. Assuming convex and Lipschitz continuous cost functions and constraints, we derive a tractable reformulation of the distributionally robust SLS (DR-SLS) problem by leveraging tools from robust control and distributionally robust optimization (DRO). Numerical experiments validate the performance and robustness of the proposed approach.
\end{abstract}


\section{Introduction}
The growing complexity of modern control systems has increasingly challenged the applicability and effectiveness of traditional control strategies. In contrast, the advancement in onboard computational capabilities and numerical solvers has opened up new opportunities for real-time optimization-based control approaches. Among them, Model Predictive Control (MPC) has emerged as a powerful and versatile framework, thanks to its ability to explicitly handle multi-objective performance criteria while systematically accounting for system constraints. MPC's flexibility in incorporating prediction models, actuator limitations, and safety requirements makes it particularly well-suited for complex and safety-critical applications such as autonomous systems, robotics, and energy management \cite{rawlings2002tutorial,schwenzer2021review,vazquez2014model}. These advantages have driven a growing interest in both theoretical developments and practical implementations of MPC across various domains.

One of the fundamental challenges in control design, including MPC, is the effective handling of uncertainties. In practice, systems are inevitably subject to various sources of uncertainty, such as model mismatch and external disturbances. The performance and reliability of a control scheme critically depend on how well these uncertainties are addressed.

Broadly speaking, methodologies for dealing with uncertainty can be categorized into two main approaches. The robust approach treats uncertainties as belonging to predefined deterministic sets and seeks to optimize control performance under the worst-case scenario \cite{bertsimas2022}. While this method is computationally efficient in general and does not require explicit knowledge of the uncertainty distributions, it often leads to overly conservative solutions. Alternatively, the stochastic approach models uncertainties probabilistically and aims to optimize the expected performance \cite{birge2011introduction}. This can reduce conservatism but typically relies on strong assumptions about the underlying uncertainty distributions—closed-form or tractable solutions are generally available only in specific cases. To address this limitation, the randomized (sample-based) approach has gained traction when a large number of uncertainty samples are available. It approximates the stochastic optimization problem by solving a large-scale optimization problem constructed from sampled scenarios \cite{campi2018introduction,calafiore2006scenario}.

To leverage the complementary strengths of both robust and stochastic control methods, the distributionally robust (DR) framework has been proposed. In this setting, the uncertainty is assumed to follow an unknown probability distribution that lies within a specified ambiguity set, which captures plausible distributions based on prior knowledge or data. The control objective is then to optimize the expected performance under the worst-case distribution within this ambiguity set \cite{mohajerin2018data,rahimian2019distributionally,kuhn2025distributionally}. Recent studies, e.g., \cite{van2015distributionally, li2021distributionally, taskesen2023distributionally, aolaritei2023wasserstein, mcallister2024distributionally, hakobyan2024wasserstein}, have explored the integration of DR methods into predictive control schemes. These works primarily focus on mitigating the impact of additive disturbances by assuming perfect knowledge of the system dynamics and optimizing open-loop control actions accordingly. 

To enhance robustness and reduce conservatism in optimal control design, closed-loop feedback policies are generally preferred over open-loop control actions in both robust optimal decision-making and Model Predictive Control (MPC) frameworks \cite{bertsimas2022, borrelli2003constrained}. However, incorporating feedback policies directly into online MPC optimization introduces significant computational challenges, as it results in nonconvex optimization problems--even when the system dynamics and control laws are linear. For linear systems, the system level synthesis (SLS) framework provides a parameterization of the closed-loop system. Other than the conventional MPC design focusing on optimizing controllers,  SLS-based design alternatively focuses on directly optimizing the whole closed-loop performance by solving convex optimization problems and reconstructing control policies from the optimized closed-loop response via linear transformations \cite{anderson2019system}. 
This shift in perspective has facilitated progress in distributed and robust predictive control design \cite{anderson2019system, sieber2021system, chen2020robust}.

As noted in \cite{schuepp2025system}, a key limitation of existing SLS methods is their reliance on purely linear control policies, primarily due to their tractability. However, linear policies lack the expressiveness to approximate piecewise-affine solutions, which frequently arise in Model Predictive Control (MPC) problems \cite{borrelli2003constrained}. Moreover, current linear-policy-based SLS frameworks cannot be directly extended to the affine case without incurring computational issues, see \cite{anderson2019system,furierilearning,furieri2019input,cescon2025data}. Notable exceptions are the works in \cite{schuepp2025system,furieri2022near}. In \cite{schuepp2025system}, the authors developed an SLS framework using state feedback control law under the assumption of full state measurability and perfect knowledge of system models. In \cite{furieri2022near}, a linear-quadratic SLS approach with linear output feedback was proposed, with a discussion on its extension to affine policies. However, this method is restricted only to LQR problems and to systems with additive input disturbances.

Motivated by the above limitations, this paper investigates SLS design with output feedback affine control policies for uncertain linear systems subject to model mismatch and additive disturbances. This setting presents several major challenges: 1) the use of affine output feedback policies complicates the construction of an SLS parameterization while maintaining computational tractability; 2) model mismatch and stochastic disturbances lead to inaccurate predictions of system dynamics, increasing the difficulty of designing effective feedback policies; 3) model mismatch induces distributional shifts between the true closed-loop behavior and the predicted closed-loop response, further complicating the solution of the DR control formulation.

In order to address the above challenges, we propose a novel DR-SLS method with output feedback affine policies. The contributions of this work are summarized as follows:
\begin{itemize}
    \item A novel SLS parameterization for a class of linear systems using output feedback affine control policy is proposed. Compared to existing SLS methods, the proposed scheme introduces fewer constraints and decision variables and accommodates a broader class of systems.
    \item The impact of model mismatch on the proposed SLS parameterization is analyzed. Besides, a new DR-SLS formulation is proposed to account for the influence of model mismatch and stochastic disturbances.
    \item The distribution shift between the prediction and true closed-loop responses is analyzed, highlighting its dependence on model uncertainties and the chosen SLS parameterization. An upper bound on this shift is derived in the sense of the Wasserstein metric.
    \item A computationally tractable relaxation of the proposed DR-SLS problem is developed by incorporating tools from robust SLS and DRO, resulting in a convex reformulation. Simulation results are provided to validate the effectiveness of the proposed approach.

\end{itemize}

The remainder of this paper is organized as follows. Section \ref{sec:problem_formu} introduces the control problem under consideration. Section \ref{sec:sls_output_affine} proposes an SLS design with output feedback affine policies and extends it to the case of inexact models. In Section \ref{sec:dro_sls}, a distributionally robust SLS (DR-SLS) formulation is introduced, along with a computationally tractable solution. Section \ref{sec:simu} provides simulation results that validate the effectiveness of the proposed approach. Finally, Section \ref{sec:conclusion} concludes this paper.

\textbf{Notation}: Lowercase letters denote vectors or scalars; uppercase letters denote matrices. Boldface lowercase letters represent the stacked version of the corresponding vectors. Norm $||\cdot||$ refers to the 1-norm for vectors and the induced 1-norm for matrices. The notation $A^{-1}$ denotes the inverse of matrix $A$. The operator $\text{col}(x_1,x_2,\cdots,x_n)$ represents the stacked vector $[x_1^{\rmt},\cdots,x_n^{\rmt}]^{\rmt}$. Calligraphic capital letters represent sets or collections. $\text{blkdiag}(M_1,\cdots,M_n)$ denotes the block-diagonal matrix with $(M_1,\cdots,M_n)$ on its diagonal block entries. $I_n$ represents a $n\times n$ identity matrix, and the subscript $n$ will be omitted when clear from the context.
\section{Problem Formulation}\label{sec:problem_formu}
Consider the following uncertain linear system in the so-called innovation form
\begin{subequations}\label{eq:sys_dyn}
    \begin{align}
        x_{t+1}& = Ax_t + Bu_t + Le_t,\\
        y_t & = Cx_t + Du_t + e_t,
    \end{align}
\end{subequations}
where $x_t\in\mathbb{R}^n$ is the system state, $u_t\in\mathbb{R}^m$ is the control input, $e_t\in\mathbb{R}^q$ is the innovation process, and subscript $t$ is time index. The matrices $A$, $B$, $L$, $C$ and $D$ are the respective system, input, Kalman, output and feedforward matrices, respectively. More details about the innovation form of linear systems can be found in \cite{breschi2023data,wang2025data,van2013closed}.

For this linear system, the following finite-horizon optimal control problem will be investigated:
\begin{subequations}\label{eq:opt}
\begin{align}
    \min_{\pi_t} & \sum_{t=0}^T h_t(y_t,u_t)\\
    \text{s.t. } & \text{system dynamics } \eqref{eq:sys_dyn},\\
    &g_t(y_t,u_t) \leq 0,\ t = 0,\cdots,T,\\
    &u_t = \pi_t(y_0,\cdots,y_{t-1}),
\end{align}
\end{subequations}
where $T$ is the length of prediction horizon, $h_t(\cdot)$ and $g_t(\cdot)$ are stage cost and constraint functions, respectively, $\pi(\cdot)$ is the control policy, which is restricted to be strictly causal. In our design, $h_t$ and $g_t$ are assumed to be convex and Lipschitz continuous functions.

In this work, only the output signal $y_t$ is assumed to be measurable for control design. To balance computational tractability and optimality, the following strictly causal output feedback affine control policy is considered
\begin{equation}\label{eq:control_policy}
    \pi_t(y_0,\cdots,y_{t-1}) = \sum_{k=0}^{t-1}K_{t,k}y_k + p_t,
\end{equation}
where $(K_{t,k},p_t)$ are control policy parameters to be optimized. To ensure the causality of the control policy, the control input $u_t$ is only parameterized via the output signals $y_k\ (k<t)$ prior to the time instant $t$.

It should be mentioned that directly solving the optimization problem \eqref{eq:opt} and considering the control policy \eqref{eq:control_policy} is computationally challenging since it leads to a nonconvex optimization problem w.r.t. the control policy parameters $K_{t,k}$ and $p_t$. Instead of directly optimizing the control parameters, the SLS parameterization provides an alternative solution by optimizing the closed-loop responses via solving a convex optimization problem. Then, the control parameters are computed based on the linear relationship between closed-loop responses and control policy parameters. 

In the upcoming sections, the SLS design for system \eqref{eq:sys_dyn} under control law \eqref{eq:control_policy} will be investigated. It is important to note that incorporating the affine control policy introduces computational challenges for existing SLS methods. Specifically, the presence of the bias term in the affine policy renders a direct extension of the existing SLS methods developed for linear control policies nonconvex, see \cite{furierilearning,furieri2019input}, and thus computationally demanding. In addition to developing an SLS design for system \eqref{eq:sys_dyn} under the assumption of perfect knowledge of system models, we further investigate the impact of model uncertainties for the proposed SLS design. To address the influence of model uncertainties, we propose a novel DR-SLS framework to improve the robustness of the optimal control solution against model mismatch and stochastic disturbances.

\section{System Level Synthesis With Output Feedback Affine Control Policy}\label{sec:sls_output_affine}

In this section, the SLS-based method for solving \eqref{eq:opt} will be presented. Firstly, a new SLS parameterization for the linear system with output feedback affine control policy is derived. Then, the effects of model mismatch on the proposed SLS formulation are analyzed.

\subsection{System Level Synthesis With Exact Model}
System \eqref{eq:sys_dyn} can be rewritten as the following predictor form\begin{subequations}\begin{align}
    x_{t+1}& = (A - LC)x_t + (B - LD) u_t + Ly_t,\\
    e_t & = y_t -Cx_t,
\end{align}
\end{subequations}
and for any $\tau \in \mathbb{Z}_+$ (also called past horizon), we see that
\begin{equation}\label{eq:init_state_estimate}
    x_t = (A-LC)^{\tau}x_{t-\tau} + \sum_{k=1}^{\tau}\left[\Psi_k^{u}u(t-k) + \Psi_k^{y}y(t-k)\right],
\end{equation}
where $\Psi_k^u = (A-LC)^{k-1}(B - LD)$ and $\Psi_k^y = (A-LC)^{k-1}L$.

Supposing that $(A-LC)$ is Schur stable and denoting $\lambda_{\max}$ ($|\lambda_{\max}| < 1$) as the eigenvalue of $(A-LC)$ with the largest absolute value, it was shown in \cite{breschi2023data} that
\begin{equation}
    x_t = \Psi
    \begin{bmatrix}
        \mathbf{u}_t^{-}\\
        \mathbf{y}_t^{-}
    \end{bmatrix}
    + \underbrace{O(|\lambda_{\max}|^{\tau})}_{\rightarrow 0 \text{ as } \tau\rightarrow +\infty}
\end{equation}
where $\Psi:= [\Psi_{\tau}^u,\cdots,\Psi_{1}^u,\Psi_{\tau}^{y},\cdots, \Psi_{1}^y]$, $\mathbf{u}_t^{-} = \text{col}(u_{t-\tau},\cdots,u_{t-1})$, and $\mathbf{y}_t^{-} = \text{col}(y_{t-\tau},\cdots,y_{t-1})$.

By selecting $\tau$ sufficiently large such that $O(|\lambda_{\max}|^{\tau})$ can be neglected, and denoting the current time instant as $t = 0$, we have
\begin{subequations}\label{eq:7}
    \begin{align}
        \mathbf{x} & = \mathbf{T}_u\mathbf{u} + \mathbf{T}_e\mathbf{e} + \Gamma \Psi 
        \begin{bmatrix}
            \mathbf{u}_0^{-}\\
            \mathbf{y}_0^{-}
        \end{bmatrix}
        , \label{eq:x_org} \\
    \mathbf{y} & = \mathbf{C}\mathbf{x} + \mathbf{D}\mathbf{u} + \mathbf{e},\label{eq:y_org}
    \end{align}
\end{subequations}
where $\mathbf{x} = \text{col}(x_0,x_1,\cdots,x_T)$ is the stacked state vector over the prediction horizon, $\mathbf{y} := \text{col}(y_0,y_1,\cdots,y_T)$ the stacked output vector, $\mathbf{u}=\text{col}(u_0,u_1,\cdots,u_T)$ the stacked control input, $\mathbf{e}:= \text{col}(e_0,\cdots,e_{T})$ the stacked innovation process, the matrices $\mathbf{C}$, $\Gamma$ and the Toeplitz matrices $\mathbf{T}_u$ and $\mathbf{T}_e$ are defined as
\begin{subequations}
    \begin{align}
\mathbf{C} &= \text{blkdiag}(\underbrace{C,\cdots,C}_{T+1 \text{ times}}),\ 
\mathbf{D} = \text{blkdiag}(\underbrace{D,\cdots,D}_{T+1 \text{ times}}),\\
    \Gamma &= 
    \begin{bmatrix}
    I\\
        A\\
        \vdots\\
        A^{T}
    \end{bmatrix},\ 
\mathbf{T}_u =
    \begin{bmatrix}
        0 & 0 & 0 & \cdots & 0\\
        B & 0 & 0 & \cdots & 0\\
        \vdots & \vdots & \ddots & \ddots & 0\\
        A^{T-1}B & A^{T-2}B & \cdots & B & 0
    \end{bmatrix},\\
    \mathbf{T}_e & =
    \begin{bmatrix}
        0 & 0 & 0 & \cdots & 0\\
        L & 0 & 0 & \cdots & 0\\
        \vdots & \vdots & \ddots & \ddots & 0\\
        A^{T-1}L & A^{T-2}L & \cdots & L & 0
    \end{bmatrix}.
    \end{align}
\end{subequations}
\noindent Combining \eqref{eq:x_org} and \eqref{eq:y_org} gives 
\begin{equation}\label{eq:output_stacked}
    \mathbf{y} = \mathbf{Gu} + \mathbf{y}_0 + \Theta\mathbf{e},
\end{equation}
where 
\begin{subequations}\label{eq:G_def}
\begin{align}
\mathbf{G} =& \mathbf{CT}_u + \mathbf{D},\\ 
\mathbf{y}_0 =& \mathbf{C}\Gamma\Psi\begin{bmatrix}
    \mathbf{u}_0^{-}\\
    \mathbf{y}_0^{-}
\end{bmatrix},\\ 
\Theta =& \mathbf{CT}_e + I.
\end{align}
\end{subequations}
It is worth noting that since $\mathbf{T}_u$ is block strictly-lower-triangular, and $\mathbf{C}$ and $\mathbf{D}$ are block diagonal matrices, the matrix $\mathbf{G}$ is a block lower-triangular matrix.

\remarknew{By considering the system dynamics in the innovation form, our proposed scheme applies to a broader class of systems compared to the approach in \cite{furieri2022near,furierilearning}, which is limited to the case of dynamic disturbances modeled as input disturbances. Additionally, the use of \eqref{eq:init_state_estimate} eliminates the need of knowledge of the initial state when predicting system outputs over the prediction horizon. Instead, only past input-output data over a finite horizon $\tau$ are required. A common heuristic for selecting $\tau$ is to choose it large enough such that the term $(A - LC)^{\tau}x_{t - \tau}$ becomes negligible, which is standard in many subspace algorithms\cite{breschi2023data,chiuso2007role}.
}

According to the control policy in \eqref{eq:control_policy}, the control inputs over the prediction horizon satisfy
\begin{equation}\label{eq:input_stacked}
    \mathbf{u} = \mathbf{Ky } + \mathbf{p},
\end{equation}
where $\mathbf{K}$ is defined as
\begin{equation}\label{eq:K_mat}
    \mathbf{K} =
    \begin{bmatrix}
        0 & 0 & \cdots & 0\\
        K_{1,0} & 0 & \cdots & 0\\
        \vdots & \ddots & \ddots & 0\\
        K_{T,0} & \cdots & K_{T,T-1} & 0
    \end{bmatrix},
\end{equation}
which is block strictly-lower-triangular, and $\mathbf{p}:=\text{col}(p_0,\cdots,p_T)$ is the bias term vector.

Combing \eqref{eq:output_stacked} and \eqref{eq:input_stacked} gives the following closed-loop representation
\begin{eqnarray}\label{eq:closed-loop_sls}
    \begin{bmatrix}
        \mathbf{y}\\
        \mathbf{u}
    \end{bmatrix} = 
    \begin{bmatrix}
        (I - \mathbf{GK})^{-1} & (I - \mathbf{GK})^{-1}\mathbf{Gp}\\
        \mathbf{K}(I - \mathbf{GK})^{-1} & (I - \mathbf{KG})^{-1}\mathbf{p}
    \end{bmatrix}
    \begin{bmatrix}
        \mathbf{y}_0 + \Theta\mathbf{e}\\
        1
    \end{bmatrix}.
\end{eqnarray}
Since $\mathbf{G}$ is block lower-triangular and $\mathbf{K}$ is block strictly-lower-triangular, it can be guaranteed that $(I - \mathbf{GK})$ and $(I - \mathbf{KG})$ are always invertible.

Based on the above analysis, we have the following theorem.

\textbf{Theorem} 1: Consider the linear system \eqref{eq:sys_dyn} and the affine control policy \eqref{eq:control_policy}, the following statements hold
\begin{enumerate}
\item the affine subspace defined by 
\begin{equation}\label{eq:sls}
    \begin{bmatrix}
        I & -\mathbf{G}
    \end{bmatrix}
    \begin{bmatrix}
        \Phi_y & \phi_y\\
        \Phi_u & \phi_u
    \end{bmatrix}=
    \begin{bmatrix}
        I & 0
    \end{bmatrix}
\end{equation}
with $\Phi_y$ and $\Phi_u$ satisfying the following block-lower-triangular and block-strictly-lower-triangular structures, respectively,
\begin{subequations}\label{eq:phi_block_lower}
    \begin{align}
    \Phi_y = &
    \begin{bmatrix}
        \Phi_y^{0,0} & 0 & \cdots & 0\\
        \Phi_y^{1,0} & \Phi_y^{1,1} & \cdots & 0\\
        \vdots & \ddots & \ddots & 0\\
        \Phi_y^{T,0} & \Phi_y^{T,1} & \cdots & \Phi_y^{T,T}
    \end{bmatrix}\\
        \Phi_u = &
    \begin{bmatrix}
        0 & 0 & \cdots & 0\\
        \Phi_u^{1,0} & 0 & \cdots & 0\\
        \vdots & \ddots & \ddots & 0\\
        \Phi_u^{T,0} & \cdots & \Phi_u^{T,T-1} & 0 
    \end{bmatrix}
    \end{align}
\end{subequations}
parameterizes all possible closed-loop output-input responses as
\begin{equation}\label{eq:sls_param}
    \begin{bmatrix}
        \mathbf{y}\\
        \mathbf{u}
    \end{bmatrix}=
    \begin{bmatrix}
        \Phi_y & \phi_y\\
        \Phi_u & \phi_u
    \end{bmatrix}
    \begin{bmatrix}
        \mathbf{y}_0 + \Theta\mathbf{e}\\
        1
    \end{bmatrix}.
\end{equation} Namely, all closed-loop response parameterizations in the format of \eqref{eq:sls_param} with $\Phi_y$ and $\Phi_u$ having structures shown in \eqref{eq:phi_block_lower} satisfy \eqref{eq:sls}. 
\item for any matrices $\{\Phi_y,\Phi_u\}$ and vectors $\{\phi_y,\phi_u\}$ satisfying \eqref{eq:sls} and \eqref{eq:phi_block_lower}, the corresponding control law $\mathbf{u}=\mathbf{Ky+p}$ with
\begin{subequations}\label{eq:control_law}
    \begin{align}
    \mathbf{K} =& \Phi_u\Phi_y^{-1},\\
    \mathbf{p} =& \phi_u - \mathbf{K}\phi_y
    \end{align}
\end{subequations}
achieves the desired closed-loop output-input responses in \eqref{eq:sls_param}.
\end{enumerate}

\textbf{Proof}: Proof of 1): It can be readily verified that the closed-loop responses shown in \eqref{eq:closed-loop_sls} with any strictly-lower-block-triangular matrix $\mathbf{K}$ satisfy \eqref{eq:sls} and \eqref{eq:phi_block_lower}.

Proof of 2): Firstly, notice that the linear constraint in \eqref{eq:sls} imposes that all diagonal blocks of $\Phi_y$ are identity matrices, which implies that $\Phi_y$ is invertible. Then, to complete the proof, it is sufficient to show that with the control law \eqref{eq:control_law}, the following equalities hold
\begin{subequations}\label{eq:sls_relationships}
    \begin{align}
    (I -\mathbf{GK})^{-1} = \Phi_y,&\ \mathbf{K}(I - \mathbf{GK})^{-1} = \Phi_u, \label{eq:cond_1}\\
    (I - \mathbf{GK})^{-1}\mathbf{Gp} = \phi_y, & \ (I - \mathbf{KG})^{-1}\mathbf{p} = \phi_u.\label{eq:cond_2}
    \end{align}
\end{subequations}
The proof of the two equalities in \eqref{eq:cond_1} follows a standard SLS design with linear policies, as shown in \cite{furieri2022near,furierilearning}. 
By substituting $\mathbf{K} = \Phi_u\Phi_y^{-1}$, we obtain
\begin{subequations}\label{eq:sls_term1_proof}
    \begin{align}
    (I - \mathbf{GK})^{-1} =& (I - \mathbf{G}\Phi_u\Phi_y^{-1})^{-1}\\
    =& (I - (\Phi_y - I)\Phi_y^{-1})^{-1} = \Phi_y,
    \end{align}
\end{subequations}
where the second equality exploits the property $\mathbf{G}\Phi_u = \Phi_y - I$ in \eqref{eq:sls}. Similarly, one has
\begin{subequations}
\begin{align*}
    \mathbf{K}(I - \mathbf{GK})^{-1} = &
    \Phi_u\Phi_y^{-1} \Phi_y = \Phi_u.
\end{align*}
\end{subequations}

Unlike in the existing SLS design with linear policies, the parameters $\phi_y$ and $\phi_u$ in \eqref{eq:cond_2} quantify the contribution of the bias term of the affine policy to the closed-loop output and input responses, respectively. Substituting $\mathbf{p} = \phi_u - \mathbf{K}\phi_y$ into the first equality in \eqref{eq:cond_2} gives
\begin{subequations}\label{eq:phi_y_res}
\begin{align}\hspace{-5pt}
(I - \mathbf{GK})^{-1}\mathbf{Gp} =& (I - \mathbf{GK})^{-1}\mathbf{G}(\phi_u - \mathbf{K}\phi_y)\\
= & (I - \mathbf{GK})^{-1}(\phi_y - \mathbf{GK}\phi_y) = \phi_y,
\end{align}
\end{subequations}
where the second equality builds upon the property $\mathbf{G}\phi_u = \phi_y$ in \eqref{eq:sls}.

Similarly, it can be derived that
\begin{subequations}
    \begin{align*}
    (I - \mathbf{KG})^{-1}\mathbf{p} =& (I + \mathbf{K}(I - \mathbf{GK})^{-1}\mathbf{G})\mathbf{p}\\
    =&\mathbf{p} + \mathbf{K}(I - \mathbf{GK})^{-1}\mathbf{Gp} \\
    = & \mathbf{p} + \mathbf{K}\phi_y = \phi_u,
    \end{align*}
\end{subequations}
where the first equality exploits the matrix inverse lemma (Sherman-Morrison-Woodbury formula) \cite{hager1989updating}, the second equality adopts the results derived in \eqref{eq:phi_y_res}, and the last equality utilizes the relationship in \eqref{eq:control_law}. This completes the proof.\hfill$\square$

\remarknew{Theorem 1 leads to a concise SLS design with output feedback affine control policy. In the existing literature, almost all SLS designs only consider linear feedback policies, based either on state or output signals, since linear policies are less theoretically challenging to deal with in the SLS framework. Compared with linear policies, affine policies show increased expressiveness and can lead to improved control performance. The exception considering the SLS design with output feedback affine policy is the work proposed in \cite{furieri2022near}. It is worth noting that, in comparison with \cite{furieri2022near}, the SLS design proposed in Theorem 1 provides a more concise design with fewer constraints and decision variables, and avoids computing an extra matrix inverse. 
}

\revision{
\remarknew{The existing SLS design approach with linear policies, e.g., \cite{furieri2022near}, could be extended to affine policies by augmenting the output signal and feedback matrix as $\bar{y}_t:=[y_t^\rmt, 1]^\rmt$ and $\bar{K}_{t,k} = [K_{t,k},p_k^\rmt]$, respectively. However, this extension introduces potential issues. Specifically, considering the augmented output signal and feedback gains causes the matrix $\mathbf{K}$ in \eqref{eq:K_mat} to be block-lower-triangular rather than strictly block-lower-triangular. Consequently, the invertibility of $(I-\mathbf{GK})$, which is required in \eqref{eq:sls_relationships}, is no longer guaranteed. However, if the system \eqref{eq:sys_dyn} is strictly proper, i.e., $D = 0$, the matrix $\mathbf{G}$ defined in \eqref{eq:G_def} becomes strictly-block-lower-triangular, and the existing approach with linear control policy is applicable to the extended system.
}}

For notational brevity, we denote the closed-loop output-input vector as $\bm{\eta}:=\text{col}(\mathbf{y},\mathbf{u})$ in the remaining parts of this work. Based on the SLS formulation in Theorem 1, the stochastic version of the finite-horizon predictive control problem \eqref{eq:opt} can be formulated as
\begin{subequations}\label{eq:sls_opt_nominal}
    \begin{align}
    \min_{\Phi_y,\Phi_u,\phi_y,\phi_u}\ & \mathbb{E}^{\mathbf{\mathbf{\bm{\eta}}}\sim\mathbb{P}^{\mathcal{M}}}\left[h(\bm{\eta})\right] \\
    \text{s.t. }& \text{constraints }\eqref{eq:sls} \text{ and } \eqref{eq:phi_block_lower},\\
    & \mathbb{E}^{\mathbf{\bm{\eta}}\sim\mathbb{P}^{\mathcal{M}}}\left[g(\bm{\eta})\right] \leq 0,\\
    & \bm{\eta} = \begin{bmatrix}
        \Phi_y & \phi_y\\
        \Phi_u & \phi_u
    \end{bmatrix}
    \begin{bmatrix}
        \mathbf{y}_0 + \Theta\mathbf{e}\\
        1
    \end{bmatrix},
    \end{align}
\end{subequations}
where $\mathbb{P}^{\mathcal{M}}$ denotes the probability distribution of the closed-loop output-input response $\bm{\eta}$ with the system model $\mathcal{M}$, which influences the parameters $(\mathbf{G},\Theta,\mathbf{y}_0,\mathbf{e})$, the functions $h(\cdot)$ and $g(\cdot)$ are constructed based on $h_t(\cdot)$ and $g_t(\cdot)$ defined in \eqref{eq:opt}, respectively. The stochastic uncertainty considered in the above formulation is the innovation process vector $\mathbf{e}$.

Suppose that the system model $\mathcal{M}$ is perfectly available and the stochastic uncertainty $\mathbf{e}$ can be replaced by its nominal prediction. By exploiting the SLS-based formulation proposed in Theorem 1, the deterministic version of the above optimization problem leads to a convex optimization problem, and the corresponding control policy can be computed via \eqref{eq:control_law}.

\subsection{SLS Parameterization with Model Mismatch}\label{sec:sls_uncertainty}
The SLS formulation proposed in Section \ref{sec:sls_output_affine}.A is based on the assumption that the system model, i.e., the set of parameters $(\mathbf{G},\Theta,\mathbf{y}_0,\mathbf{e})$, is perfectly known. However, obtaining the exact model information is nontrivial and challenging. It is more practical to consider that only an approximate model is available for SLS design. 

It is worth noting that among the above-listed parameters, $\mathbf{G}$ plays the most important role in influencing the SLS design since it will directly influence the SLS design via \eqref{eq:sls}. In this section, the SLS formulation proposed in Theorem 1 of Section \ref{sec:sls_output_affine}.A is extended to the case where only the approximated parameter $\hat{\mathbf{G}}$ is available. For ease of notation, the approximation error of $\mathbf{G}$ is defined as $\Delta := \mathbf{G} - \hat{\mathbf{G}}$.


\textbf{Theorem} 2: Assuming a nominal approximation parameter $\hat{\mathbf{G}}$ is adopted for SLS-based design in \eqref{eq:sls} and \eqref{eq:control_law}, the true closed-loop input-output response $\bm{\eta}: = \text{col}(\bm{y},\bm{u})$ with the true parameters $(\mathbf{G},\mathbf{y}_0,\Theta,\mathbf{e})$ is
\begin{equation}\label{eq:sls_true}
    \bm{\eta} = \Phi(\mathbf{y}_0 + \Theta\mathbf{e}) + \phi
\end{equation}
with
\begin{subequations}\label{eq:phi_def}
    \begin{align}
    \Phi := 
    \begin{bmatrix}
        \Phi_y\\
        \Phi_u
    \end{bmatrix} =& 
    \begin{bmatrix}
        \hat{\Phi}_y\\
        \hat{\Phi}_u
    \end{bmatrix}\underbrace{(I - \Delta \hat{\Phi}_u)^{-1}}_{R_{\Phi}},\label{eq:def_Phi}\\
    \phi := 
    \begin{bmatrix}
        \phi_y\\
        \phi_u
    \end{bmatrix} =& \underbrace{
    \begin{bmatrix}
        I & \hat{\Phi}_y(I - \Delta\hat{\Phi}_u)^{-1}\Delta\\
        0 & I + \hat{\Phi}_u(I - \Delta\hat{\Phi}_u)^{-1}\Delta
    \end{bmatrix}}_{R_{\phi}}
    \begin{bmatrix}
        \hat{\phi}_y\\
        \hat{\phi}_u
    \end{bmatrix},\label{eq:def_phi}
    \end{align}
\end{subequations}
where $\{\Phi_y,\Phi_u,\phi_y,\phi_u\}$ is the true closed-loop response parametrization, $\{\hat{\Phi}_y,\hat{\Phi}_u,\hat{\phi}_y,\hat{\phi}_u\}$ is the predictive closed-loop response parameterization derived via \eqref{eq:sls} by replacing $\mathbf{G}$ with $\hat{\mathbf{G}}$, and $(R_{\Phi}, R_{\phi})$ are the uncertainties of closed-loop parameterization introduced by the model mismatch $\Delta$.

\textbf{Proof}: Based on Theorem 1, we have 
\begin{equation*}
    \begin{bmatrix}
        I & -\mathbf{\hat{G}} 
    \end{bmatrix}
    \begin{bmatrix}
        \hat{\Phi}_y \\
        \hat{\Phi}_u
    \end{bmatrix} = I,
\end{equation*}
which further leads to
\begin{equation*}
    \begin{bmatrix}
        I & -\mathbf{G}
    \end{bmatrix}
    \begin{bmatrix}
        \hat{\Phi}_y \\
        \hat{\Phi}_u
    \end{bmatrix}
     = I - \Delta\hat{\Phi}_u.
\end{equation*}
Consequently, we arrive at
\begin{equation}\label{eq:true_sls_cond1}
    \begin{bmatrix}
        I & -\mathbf{G}
    \end{bmatrix}
    \begin{bmatrix}
        \hat{\Phi}_y\\
        \hat{\Phi}_u
    \end{bmatrix}\left(I - \Delta \hat{\Phi}_u\right)^{-1} = I,
\end{equation}
where the invertibility of $(I - \Delta\hat{\Phi}_u)$ is guaranteed since $\Delta$ is block-lower-triangular and $\hat{\Phi}_u$ is strictly-block-lower-triangular. The above equality provides an option for selecting $\Phi_y = \hat{\Phi}_y(I - \Delta\hat{\Phi}_u)^{-1}$ and $\Phi_u = \hat{\Phi}_u(I - \Delta\hat{\Phi}_u)^{-1}$ that ensures $\mathbf{K} = \Phi_u\Phi_y^{-1} = \hat{\Phi}_u\hat{\Phi}_y^{-1}$.

Since $\hat{\phi}_y$ and $\hat{\phi}_u$ satisfy
\begin{equation}
\begin{bmatrix}
    I & -\hat{\mathbf{G}}
\end{bmatrix}
\begin{bmatrix}
    \hat{\phi}_y\\
    \hat{\phi}_u
\end{bmatrix} = 0,
\end{equation}
it follows that
\begin{equation}\label{eq:phi_1}
    \begin{bmatrix}
        I & -\mathbf{G}
    \end{bmatrix}
    \begin{bmatrix}
        \hat{\phi}_y\\
        \hat{\phi}_u
    \end{bmatrix} = 
    -\Delta\hat{\phi}_u.
\end{equation}

Since $(\Phi_y,\Phi_u,\phi_y,\phi_u)$ are defined as true closed-loop parameterizations under the control policy parameters $(\mathbf{K},\mathbf{p})$ derived via \eqref{eq:control_law}, Theorem 1 implies the following equalities
\begin{subequations}\label{eq:phi_2}
\begin{align}
    & \begin{bmatrix}
        I & -\mathbf{G}
    \end{bmatrix}
 \begin{bmatrix}
     \phi_y\\
     \phi_u
 \end{bmatrix} = 0,\\
    & \mathbf{p} = \phi_u - \mathbf{K}\phi_y = 
\hat{\phi}_u - \mathbf{K}\hat{\phi}_y,\label{eq:p_equal}
\end{align}
\end{subequations}
where the relationship in \eqref{eq:p_equal} is built upon the fact that the control policy parameter $\mathbf{p}$ is identical in both cases of the true closed-loop parameterization and nominal closed-loop parameterization.

Denoting $r_y: = \phi_y - \hat{\phi}_y$ and $r_u := \phi_u - \hat{\phi}_u$, and combining \eqref{eq:phi_1} and \eqref{eq:phi_2} results in
\begin{equation}\label{eq:phi_cond}
\begin{cases}
    r_u - \mathbf{K}r_y = 0\\
    r_y - \mathbf{G}r_u = \Delta\hat{\phi}_u
    \end{cases}
\end{equation}
Solving \eqref{eq:phi_cond} gives
\begin{subequations}\label{eq:true_sls_cond2}
    \begin{align}
    r_u =\ &\mathbf{K}(I - \mathbf{GK})^{-1}\Delta\hat{\phi}_u = \Phi_u\Delta \hat{\phi}_u,\\
    r_y =\ & (I - \mathbf{GK})^{-1}\Delta \hat{\phi}_u = \Phi_y\Delta\hat{\phi}_u,
    \end{align}
\end{subequations}
where the relationships in \eqref{eq:sls_relationships} are exploited.

Combining the results in \eqref{eq:true_sls_cond1} and \eqref{eq:true_sls_cond2} and considering the definition of $(r_y,r_u)$ leads to the closed-loop parameterizations in \eqref{eq:sls_true} and \eqref{eq:phi_def}. For these parameters, it can be verified that the condition \eqref{eq:sls} in Theorem 1 is satisfied, which justifies that the parameters in \eqref{eq:phi_def} are true closed-loop parameterizations. 

Furthermore, based on the analysis shown in \eqref{eq:true_sls_cond1}-\eqref{eq:true_sls_cond2}, it can be verified that the parameters $\{\Phi_y,\Phi_u,\phi_y,\phi_u\}$ together with the true model $\mathbf{G}$ satisfy \eqref{eq:sls}, and 
\begin{subequations}
    \begin{align}
        \mathbf{K} =\ & \Phi_u \Phi_y^{-1} = \hat{\Phi}_u\hat{\Phi}_y^{-1},\\
        \mathbf{p} =\ & \phi_u - \mathbf{K}\phi_y = \hat{\phi}_u - \mathbf{K}\hat{\phi}_y.
    \end{align}
\end{subequations}

Given the results in Theorem 1, we can conclude that $\{\Phi_y,\Phi_u,\phi_y,\phi_u\}$ defined in \eqref{eq:phi_def} are the true closed-loop response parameterizations under the control policy $\mathbf{u} = \mathbf{Ky+p}$, where the control policy parameters $(\mathbf{K},\mathbf{p})$ are derived based the SLS design with the nominal model $\hat{\mathbf{G}}$ via \eqref{eq:sls} and \eqref{eq:control_law}. This completes the proof. \hfill $\square$

\section{Distributionally Robust SLS Design and Tractable Reformulation}\label{sec:dro_sls}
Based on the results shown in Section \ref{sec:sls_output_affine}, this section proposes a novel DR-SLS framework to systematically enhance the robustness of the optimal control solution against model mismatch and additive stochastic disturbances.
\subsection{Distributionally Robust SLS Formulation}

Solving the SLS-based problem \eqref{eq:sls_opt_nominal} is generally intractable when assuming the existence of model mismatch for parameters $(\mathbf{G},\Theta,\mathbf{y}_0)$, and uncertain probability distribution $\mathbb{P}_{\mathbf{e}}$. One possible solution for solving \eqref{eq:sls_opt_nominal} is adopting the sample average approximation (SAA) and certainty equivalence (CE) principles, in which the unknown system parameters $(\mathbf{G},\Theta,\mathbf{y}_0)$ are replaced by their nominal values $(\hat{\mathbf{G}},\hat{\Theta},\hat{\mathbf{y}}_0)$ and the true innovation process distribution $\mathbb{P}_{\mathbf{e}}$ is approximated by its empirical counterpart $\overline{\mathbb{P}}_{\mathbf{e}}$. In this way, the true closed-loop response distribution $\mathbb{P}^{\mathcal{M}}$ is replaced by its empirical approximated distribution $\overline{\mathbb{P}}^{\hat{\mathcal{M}}}$.

Although this solution is straightforward and easy to implement, its efficacy is built upon the expectation that the true closed-loop distribution $\mathbb{P}^{\mathcal{M}}$ can be properly approximated by $\overline{\mathbb{P}}^{\hat{\mathcal{M}}}$. However, the complex dependency of the closed-loop system behaviour on the true model $\mathcal{M}$, as opposed to the nominal predictive model $\hat{\mathcal{M}}$, undermines the justification for applying CE and SAA. Moreover, when only a limited number of uncertainty samples are available, the SAA approach may lead to a poor out-of-sample performance.

In the following, a distributionally robust SLS (DR-SLS) design will be presented, considering the uncertainties introduced by the model mismatch between $(\mathbf{G},\Theta, \mathbf{y}_0)$ and $(\hat{\mathbf{G}},\hat{\Theta},\hat{\mathbf{y}}_0)$ as well as the estimation error of the empirical innovation process distribution $\overline{\mathbb{P}}_{\mathbf{e}}$. The approximation errors of $(\mathbf{G},\Theta,\mathbf{y}_0)$ are denoted as
$\Delta:= \mathbf{G} - \hat{\mathbf{G}}$, $\tilde{\Theta} := \Theta - \hat{\Theta}$, and $\tilde{\mathbf{y}}_0:= \mathbf{y}_0 - \hat{\mathbf{y}}_0$, respectively.
Before presenting the main results of this section, the following definitions are introduced.

\textit{Definition} 1. (Wasserstein Distance\cite{mohajerin2018data}) The Wasserstein distance $d_W: \mathcal{H}(\Xi)\times\mathcal{H}(\Xi)\rightarrow \mathbb{R}$ is defined as
\begin{equation*}
    d_W(\mathbb{Q}_1,\mathbb{Q}_2) := \inf_{\Pi} \left\{\int_{\Xi^2}||\xi_1-\xi_2||\Pi(d\xi_1,d\xi_2)\right\},
\end{equation*}
where $\Pi$ is a joint distribution of the random variables $\xi_1$ and $\xi_2$ with marginal distributions $\mathbb{Q}_1$ and $\mathbb{Q}_2$, respectively, and $\mathcal{H}(\Xi)$ is the space of all distributions $\mathbb{Q}$ supported on $\Xi$ such that $\int_{\Xi}||\xi||\mathbb{Q}(d\xi)< +\infty$.

\textit{Definition} 2. (Wasserstein Ambiguity Set \cite{mohajerin2018data}) The Wasserstein ambiguity set $\mathcal{B}_{\varepsilon}(\mathbb{Q}_0)$ is defined as the ball of radius $\varepsilon$ centered at the empirical distribution $\mathbb{Q}_0$
\begin{equation*}
    \mathcal{B}_{\varepsilon}(\mathbb{Q}_0) :=\left\{\mathbb{Q}\in\mathcal{H}(\Xi)\mid d_W(\mathbb{Q},\mathbb{Q}_0) \leq \varepsilon\right\}.
\end{equation*}

For notational brevity, the true system model and closed-loop parameterizations are compactly represented as $\mathcal{M}:=\{\mathbf{G},{\Theta},\Phi_y,\Phi_u,\phi_y,\phi_u,\mathbf{y}_0\}$. Similarly, the nominal predictive closed-loop model and the corresponding parameterizations are denoted as $\hat{\mathcal{M}}:=\{\hat{\mathbf{G}},\hat{{\Theta}},\hat{\Phi}_y,\hat{\Phi}_u,\hat{\phi}_y,\hat{\phi}_u,\hat{\mathbf{y}}_0\}$. 


\assumption{ For the innovation process $\mathbf{e}$, $N$ samples of historical realizations $\{\mathbf{e}_i\}_{i=1}^N$ are available.}

The DR-SLS problem is formulated as 
\begin{subequations}\label{eq:dro_sls}
    \begin{align}
\min_{\substack{\hat{\Phi}_y,\hat{\Phi}_u\\\hat{\phi}_y,\hat{\phi}_u}}\sup_{\mathbb{Q}\in\mathcal{B}_\varepsilon(\overline{\mathbb{P}}^{\hat{\mathcal{M}}})}\ & \mathbb{E}^{\bm{\eta}\sim\mathbb{Q}}\left[h(\bm{\eta})\right]\label{eq:dro_sls_obj}\\
    \text{s.t. } & \sup_{\mathbb{Q}\in\mathcal{B}_\varepsilon (\overline{\mathbb{P}}^{\hat{\mathcal{M}}})} \mathbb{E}^{\bm{\eta}\sim \mathbb{Q}}[g(\bm{\eta})] \leq 0,  \label{eq:dro_sls_cons}\\
    & \begin{bmatrix}
        I & -\hat{\mathbf{G}}
    \end{bmatrix}
    \begin{bmatrix}
        \hat{\Phi}_y & \hat{\phi}_y \\
        \hat{\Phi}_u & \hat{\phi}_u
    \end{bmatrix} = 
    \begin{bmatrix}
        I & 0
    \end{bmatrix}, \label{eq:nominal_sls}\\
    & \text{constraints } \eqref{eq:sls_true}  \text{ and } \eqref{eq:phi_def},
    \end{align}
\end{subequations}
where $\mathbb{Q}$ is the probability distribution of the true closed-loop output-input responses $\bm{\eta}$, $\overline{\mathbb{P}}^{\hat{\mathcal{M}}}$ is the empirical predictive probability distribution of closed-loop dynamics based on $(\hat{\mathbf{G}},\hat{\Theta}, \hat{\mathbf{y}}_0,\overline{\mathbb{P}}_{\mathbf{e}})$, $\mathcal{B}_\varepsilon(\overline{\mathbb{P}}^{\hat{\mathcal{M}}})$ is the ambiguity set of $\mathbb{Q}$.

Given $N$ samples of the innovation process $\{\mathbf{e}_i\}_{i=1}^N$, the empirical distribution $\overline{\mathbb{P}}_{\mathbf{e}}$ can be constructed as $\overline{\mathbb{P}}_{\mathbf{e}}:=\frac{1}{N}\sum_{i=1}^N\delta(\mathbf{e}_i)$, where $\delta(\cdot)$ is the Dirac distribution. Accordingly, the empirical predictive distribution of the output-input response $\bm{\eta}$ with the nominal model $\mathcal{\hat{M}}$ is
\begin{subequations}\label{eq:nominal_dis}
    \begin{align}
    \overline{\mathbb{P}}^{\hat{\mathcal{M}}} = &\ \frac{1}{N}\sum_{i=1}^N\delta(\hat{\bm{\eta}}_i), \\
    \hat{\bm{\eta}}_i =&\  \hat{\Phi}(\hat{\mathbf{y}}_0 + \hat{\Theta}\mathbf{e}_i) + \hat{\phi},
    \end{align}
\end{subequations}
where $(\hat{\Phi},\hat{\phi})$ are the predictive closed-loop response parameterization derived with $\hat{\mathbf{G}}$ via \eqref{eq:nominal_sls}, and $\hat{\bm{\eta}}_i$ is the empirical prediction of output-input response using $(\hat{\Phi},\hat{\phi},\hat{\Theta},\hat{\mathbf{y}}_0,\mathbf{e}_i)$.

\subsection{Characterization of Distribution Shift}
\revision{
In order to properly formulate and solve the DR-SLS problem in \eqref{eq:dro_sls}, the distribution shift between $\mathbb{P}^{\mathcal{M}}$ and $\overline{\mathbb{P}}^{\hat{\mathcal{M}}}$, which is measured by the Wasserstein distance $d_W(\overline{\mathbb{P}}^{\hat{\mathcal{M}}},\mathbb{P}^{\mathcal{M}})$ and the associated ambiguity set $\mathcal{B}_{\varepsilon}(\overline{\mathbb{P}}^{\hat{\mathcal{M}}})$, should be properly considered. An excessively large radius $\varepsilon$ of the ambiguity set will lead to a conservative solution, and an excessively small one would result in an aggressive and fragile solution. Importantly, the Wasserstein distance $d_W(\overline{\mathbb{P}}^{\hat{\mathcal{M}}},\mathbb{P}^{\mathcal{M}})$ is determined by the closed-loop model mismatch between $\mathcal{M}$ and $\hat{\mathcal{M}}$, and is therefore a function of the predictive closed-loop parameterization $\{\hat{\Phi}_y,\hat{\Phi}_u,\hat{\phi}_y,\hat{\phi}_u\}$. Consequently, simply setting the radius of the ambiguity set with a prescribed constant does not reflect its intrinsic dependence on the SLS parameters and can lead to unsatisfactory performance.}

Based on the results in Theorem 2, the true closed-loop output-input responses $\bm{\eta}\sim\mathbb{P}^{\mathcal{M}}$ are
\begin{equation}\label{eq:true_dis}
    \bm{\eta} = \hat{\Phi}R_{\Phi}(\mathbf{y}_0 +\Theta\mathbf{e}) + R_{\phi}\hat{\phi}.
\end{equation}
Compared with the empirical prediction in \eqref{eq:nominal_dis}, there are two factors contributing to the distribution shift: 1) the uncertainties caused by the model mismatch between $(\mathbf{G},\Theta,\mathbf{y}_0)$ and $(\hat{\mathbf{G}},\hat{\Theta},\hat{\mathbf{y}}_0)$, and 2) the approximation error between the true distribution $\mathbb{P}_{\mathbf{e}}$ and its empirical approximation $\overline{\mathbb{P}}_{\mathbf{e}}$. 

In the following, we analyze the Wasserstein metric $d_W(\overline{\mathbb{P}}_{\pi}^{\mathcal{\hat{M}}},\mathbb{P}_{\pi}^{\mathcal{M}})$, and provide an upper bound for the distribution shift. The triangle inequality of Wasserstein distance gives
\begin{equation}\label{eq:triangle}
    d_W(\overline{\mathbb{P}}^{\hat{\mathcal{M}}},\mathbb{P}^{\mathcal{M}}) \leq d_W(\overline{\mathbb{P}}^{\hat{\mathcal{M}}},\overline{\mathbb{P}}^{\mathcal{M}}) + d_W(\overline{\mathbb{P}}^{\mathcal{M}},\mathbb{P}^{\mathcal{M}}),
\end{equation}
where $\overline{\mathbb{P}}^{\mathcal{M}}$ is the empirical version of the true output-input distribution $\mathbb{P}^{\mathcal{M}}$. Namely,
\begin{subequations}\label{eq:eta_i}
    \begin{align}
        \overline{\mathbb{P}}^{\mathcal{M}} =\ & \frac{1}{N}\sum_{i=1}^N\delta(\bm{\eta}_i), \\
        \bm{\eta}_i =\ & \hat{\Phi}R_{\Phi}(\mathbf{y}_0 + \Theta \mathbf{e}_i )  + R_{\phi}\hat{\phi}.
    \end{align}
\end{subequations}
Applying \eqref{eq:triangle} decomposes the distribution shift into two parts: a) the shift caused by model mismatch $d_W(\overline{\mathbb{P}}^{\hat{\mathcal{M}}}, \overline{\mathbb{P}}^{\mathcal{M}})$, and b) the shift caused by the emprical approximation error $d_W(\overline{\mathbb{P}}^{\mathcal{M}},\mathbb{P}^{\mathcal{M}})$. The following lemmas provide upper bounds for the above distribution shifts.

\lemma{(distribution shift by model mismatch) The distribution shift caused by the model mismatch $d_W(\overline{\mathbb{P}}^{\hat{\mathcal{M}}}, \overline{\mathbb{P}}^{\mathcal{M}})$ can be upper-bounded as 
\begin{equation}
\begin{aligned}
d_W(\overline{\mathbb{P}}^{\mathcal{\hat{M}}},\overline{\mathbb{P}}^{\mathcal{M}}) \leq & \frac{1}{N} \sum_{i=1}^N\bigg(||\hat{\Phi}(R_{\Phi}-I)(\hat{\mathbf{y}}_0 + \hat{\Theta}\mathbf{e}_i)|| +\\
    & \quad  ||\hat{\Phi}R_{\Phi} (\tilde{\mathbf{y}}_0 +\tilde{\Theta}\mathbf{e}_i)||\bigg) + ||(R_{\phi} - I)\hat{\phi}||.
\end{aligned}
\end{equation}
}
\textbf{Proof}: The Wassertein distance $d_W(\overline{\mathbb{P}}^{\mathcal{\hat{M}}},\overline{\mathbb{P}}^{\mathcal{M}})$ is defined as
\begin{equation}\label{eq:shift_model_mismatch}
    d_W(\overline{\mathbb{P}}^{\mathcal{\hat{M}}},\overline{\mathbb{P}}^{\mathcal{M}}) := \inf_{\Pi}\int_{\Xi^2}||\hat{\bm{\eta}}_i - \bm{\eta}_i||\Pi(d\hat{\bm{\eta}}_i,d\bm{\eta}_i),
\end{equation}
where $\Pi$ is the joint probability distribution of $(\hat{\bm{\eta}}_i,\bm{\eta}_i)$. By selecting the suboptimal adjoint probability distribution $\Pi^*(\hat{\bm{\eta}}_i,\bm{\eta}_i) = \frac{1}{N}\delta(\hat{\bm{\eta}}_i,\bm{\eta}_i)$, and considering the expression of $\hat{\bm{\eta}}_i$ and $\bm{\eta}_i$ in \eqref{eq:nominal_dis} and \eqref{eq:eta_i} together with \eqref{eq:phi_def}, and applying the norm triangular inequality completes the proof. \hfill $\square$

\lemma{(distribution shift by empirical approximation error of $\mathbb{P}_{\mathbf{e}}$) The Wasserstein distance between the empirical closed-loop output-input distribution $\overline{\mathbb{P}}^{\mathcal{M}}$ and the real closed-loop output-input distribution $\mathbb{P}^{\mathcal{M}}$ is upper-bounded as
\begin{equation}
    d_W(\overline{\mathbb{P}}^{\mathcal{M}},\mathbb{P}^{\mathcal{M}}) \leq ||\hat{\Phi} R_{\Phi}(\hat{\Theta}+\tilde{\Theta})||d_W(\overline{\mathbb{P}}_{\mathbf{e}},\mathbb{P}_{\mathbf{e}}). 
\end{equation}
}
\textbf{Proof}: Based on the definition of Wasserstein distance, we have
\begin{small}\begin{align*}
    d_W(\overline{\mathbb{P}}^{\mathcal{M}},\mathbb{P}^{\mathcal{M}}) := &\inf_{\Pi}\int_{\Xi^2}||\bm{\eta}^{\prime} - \bm{\eta}||\Pi(d\bm{\eta^{\prime}},d{\bm{\eta}}) \\ 
    = & \inf_{\Pi}\int_{\mathcal{V}^2}||\hat{\Phi}R_{\Phi}(\hat{\Theta} + \tilde{\Theta})(\mathbf{e}^\prime - \mathbf{e})||\Pi(d\mathbf{e}^\prime,d\mathbf{e})\\
    \leq & ||\hat{\Phi}R_{\Phi}(\hat{\Theta} + \tilde{\Theta})||\underbrace{\inf_{\Pi}\int_{\mathcal{V}^2}||\mathbf{e}^{\prime} - \mathbf{e}||\Pi(d\mathbf{e}^{\prime},d\mathbf{e})}_{d_W(\overline{\mathbb{P}}_{\mathbf{e}},\mathbb{P}_{\mathbf{e}})}
\end{align*}\end{small}where the above relationships are based on the fact that the innovation process $\mathbf{e}$ is coupled with $\bm{\eta}$ only through $\Phi$ in \eqref{eq:sls_true} and \eqref{eq:phi_def}, and the sub-multiplicativity property of 1-norm. 
This completes the proof.\hfill $\square$

\proposition{The Wasserstein distance between the empirical predictive output-input distribution $\overline{\mathbb{P}}^{\hat{\mathcal{M}}}$ and the true output-input distribution $\mathbb{P}^{\mathcal{M}}$ is upper-bounded by:
\begin{small}\begin{equation}\label{eq:Wasserstein_bound}
\begin{aligned}
    d_W(\overline{\mathbb{P}}^{\hat{\mathcal{M}}},\mathbb{P}^{\mathcal{M}} ) \leq &\frac{1}{N} \sum_{i=1}^N\bigg(||\hat{\Phi}(R_{\Phi}-I)(\hat{\mathbf{y}}_0 + \hat{\Theta}\mathbf{e}_i)|| \\
    & \quad + ||\hat{\Phi}R_{\Phi}(\tilde{\mathbf{y}}_0+\tilde{\Theta}\mathbf{e}_i)||\bigg) + ||(R_{\phi} - I)\hat{\phi}|| \\
    &\quad + ||\hat{\Phi}R_{\Phi}(\hat{\Theta} + \tilde{\Theta})||d_W(\overline{\mathbb{P}}_{\mathbf{e}},\mathbb{P}_{\mathbf{e}}).
\end{aligned}
\end{equation}\end{small}
}
\textbf{Proof}: Applying the triangle inequality \eqref{eq:triangle} and the results in Lemmas 1 and 2 completes the proof. \hfill$\square$

\subsection{Tractable Reformulation}
While the DR-SLS formulation in \eqref{eq:dro_sls} could enhance the robustness of the optimal solution of the SLS design, it is challenging to solve the corresponding optimization problem when considering the Wasserstein distance \eqref{eq:Wasserstein_bound}. In this section, a tractable reformulation is proposed. To this end, the following assumptions are made.

\assumption{There exist constants $\alpha > 1$ and $\beta>0$ such that 
\begin{equation}
    \mathbb{E}^{\mathbf{e}\sim\mathbb{P}_{\mathbf{e}}}\left[e^{\beta||\mathbf{e}||^\alpha}\right] < +\infty
\end{equation}}

\assumption{The Wasserstein distance between the empirical distribution $\overline{\mathbb{P}}_{\mathbf{e}}$ and the true distribution $\mathbb{P}_{\mathbf{e}}$ is bounded as $d_W(\overline{\mathbb{P}}_{\mathbf{e}},\mathbb{P}_{\mathbf{e}}) \leq \kappa $.}

\assumption{The model mismatch is bounded as $||\Delta|| \leq \gamma_1$, $||\tilde{\Theta}|| \leq \gamma_2$, and $||\tilde{\mathbf{y}}_0|| \leq \gamma_3$.}

\assumption{The cost function $h(\bm{\eta})$ and the system constraint $g(\bm{\eta})$} in \eqref{eq:dro_sls} are convex and Lipschitz continuous, with Lipschitz constants $l_h$ and $l_g$, respectively.

\remarknew{Assumption 2 ensures that the distribution $\mathbb{P}_{\mathbf{e}}$ has light tails, which is standard in DRO problems to guarantee well-posedness. Assumption 3 imposes a bound on the Wasserstein distance. The radius $\kappa$ can be selected using, for example, cross-validation. Assumption 4 provides upper bounds on model mismatch, which is standard in robust control design (e.g., \cite{furieri2022near,furierilearning,micheli2024data}). If uncertainty bounds are originally given for the system matrices $(A,B,C,D)$, the corresponding bounds for $(\Delta, \tilde{\Theta}, \tilde{\mathbf{y}}_0)$ can be obtained via either numerical simulation or theoretical analysis as shown in \cite{liu2024stability}. Assumption 5 ensures that the cost and constraint functions are tractable by assuming convexity and Lipschitz continuity. This covers a broad class of functions, including linear functions, pointwise maximum of affine functions, and norms such as the $l_1$-norm and $l_2$-norm.

}


\textbf{Theorem} 3. Supposing Assumptions 1-5 hold, 
and letting $\gamma_1||\hat{\Phi}_u||\leq \rho$ with $\rho\in[0,1)$ and $||\hat{\Phi}_y|| \leq \sigma$ with $\sigma>0$, a convex relaxation of the DR-SLS problem in \eqref{eq:dro_sls} can be formulated as
\begin{subequations}\label{eq:dro_final}
    \begin{align}
        \min_{\overline{\varepsilon},\hat{\Phi},\hat{\phi},s_i,q_i} &\ l_h\overline{\varepsilon} + \frac{1}{N}\sum_{i=1}^Ns_i \\
        \text{s.t.} &\begin{bmatrix}
            I & -\hat{\mathbf{G}}
            \end{bmatrix}
            \begin{bmatrix}
                \hat{\Phi}_y & \hat{\phi}_y\\
                \hat{\Phi}_u & \hat{\phi}_u
            \end{bmatrix} = 
            \begin{bmatrix}
                I & 0
            \end{bmatrix},\label{eq:final_sls_1}\\
            &\hat{\bm{\eta}}_i = \begin{bmatrix}
                \hat{\Phi}_y & \hat{\phi}_y\\
                \hat{\Phi}_u & \hat{\phi}_u
            \end{bmatrix}
            \begin{bmatrix}
                \hat{\mathbf{y}}_0 + \hat{\Theta}\mathbf{e}_i\\
                1
            \end{bmatrix},\\
            & \hat{\Phi}_y \text{ and } \hat{\Phi}_u \text{ satisfy the structure in }\eqref{eq:phi_block_lower},\label{eq:final_sls_3}\\
            &h(\hat{\bm{\eta}}_i)\leq s_i,\quad \forall i\in\{0,\cdots,N\},\label{eq:final_dro_obj_1}\\
            & l_g\overline{\varepsilon}+\frac{1}{N}\sum_{i=1}^Nq_i \leq 0,\label{eq:final_dro_cons_1}\\
                &g(\hat{\bm{\eta}}_i)\leq q_i,\quad \forall i\in\{0,\cdots,N\},\label{eq:final_dro_cons_2}\\
            & \overline{\varepsilon} \geq \frac{||\hat{\Phi}||}{N(1 - \rho)}\bigg[\sum_{i=1}^N\left(\rho||\hat{\mathbf{y}}_0+\hat{\Theta}\mathbf{e}_i|| + \gamma_2||\mathbf{e}_i||\right) \notag\\ & \qquad +N\gamma_3\bigg] 
              + \left(\frac{\rho}{\gamma_1}+\sigma\right)\cdot\frac{\gamma_1}{1-\rho}||\hat{\phi}_u|| \notag\\
              &\qquad + \frac{\kappa}{1-\rho}||\hat{\Phi}||\cdot(||\hat{\Theta}||+\gamma_2),\label{eq:Wassers_distant}\\
              & ||\hat{\Phi}_y||\leq \sigma,\ \gamma_1||\hat{\Phi}_u|| \leq \rho,\ \rho\in[0,1),\label{eq:small_gain}
    \end{align}
\end{subequations}where constraints \eqref{eq:final_sls_1}-\eqref{eq:final_sls_3} are the SLS-based parameterization of closed-loop output-input response, constraint \eqref{eq:final_dro_obj_1} is for reformulating the objective function in \eqref{eq:dro_sls_obj}, and constraints \eqref{eq:final_dro_cons_1} and \eqref{eq:final_dro_cons_2} are for reformulating the constraint in \eqref{eq:dro_sls_cons}, constraint \eqref{eq:Wassers_distant} is for upper-bounding the Wasserstein distance $d_W(\overline{\mathbb{P}}^{\hat{\mathcal{M}}}, \mathbb{P}^{\mathcal{M}})$. 

\textbf{Proof}: The key step in proving Theorem 3 is to derive the computationally tractable relaxation in \eqref{eq:Wassers_distant} for the upper bound of the Wasserstein distance $d_W(\overline{\mathbb{P}}^{\hat{\mathcal{M}}},\mathbb{P}^{\mathcal{M}})$ proposed in Proposition 1. It can be observed that the right-hand side of \eqref{eq:Wasserstein_bound} consists of four terms. In the following, computationally tractable relaxations for these terms are derived.

Term 1: $||\hat{\Phi}(R_{\Phi} - I)(\hat{\mathbf{y}}_0 + \hat{\Theta}\mathbf{e}_i)||$. The norm sub-multiplicativity can be exploited to show that
\begin{equation}\label{eq:term1}
    ||\hat{\Phi}(R_{\Phi} - I)(\hat{\mathbf{y}}_0 + \hat{\Theta}\mathbf{e}_i)|| \leq ||\hat{\Phi}||\cdot||(R_{\Phi} - I)||\cdot||\hat{\mathbf{y}}_0 + \hat{\Theta}\mathbf{e}_i||.\\
\end{equation}
The constraint \eqref{eq:small_gain} together with Assumption 4 implies that $||\Delta\hat{\Phi}_u|| \leq \rho < 1$. Then, based on Neumann series formulation, norm triangle inequality, and norm sub-multiplicativity, we can derive that
\begin{subequations}\label{eq:term1_sub}
    \begin{align}
    \hspace{-6pt}||R_{\Phi} - I|| = &||\sum_{k=1}^{\infty}(\Delta\hat{\Phi}_u)^{k}||\\
    \leq& \sum_{k=1}^{\infty}||(\Delta\hat{\Phi}_u)^k||
    \leq \sum_{k=1}^{\infty}||\Delta\hat{\Phi}_u||^k
    \leq \frac{\rho}{1-\rho}.
    \end{align}
\end{subequations}
Substituting \eqref{eq:term1_sub} into \eqref{eq:term1} gives 
\begin{equation}\label{eq:term1_final}
   ||\hat{\Phi}(R_{\Phi} - I)(\hat{\mathbf{y}}_0 + \hat{\Theta}\mathbf{e}_i)|| \leq \frac{\rho}{1-\rho}||\hat{\Phi}||\cdot|| \hat{\mathbf{y}}_0+\hat{\Theta}\mathbf{e}_i||.
\end{equation}

Term 2: $||\hat{\Phi}R_{\Phi}(\tilde{\mathbf{y}}_0+\tilde{\Theta}\mathbf{e}_i)||$. Following a similar analysis as in Term 1, and exploiting Assumption 4, it can be derived that
\begin{equation}\label{eq:term2_final}
    ||\hat{\Phi}R_{\Phi}(\tilde{\mathbf{y}}_0+\tilde{\Theta}\mathbf{e}_i)|| \leq  \frac{1}{1-\rho}||\hat{\Phi}||\big(\gamma_2||\mathbf{e}_i|| + \gamma_3\big).
\end{equation}

Term 3: $||(R_{\phi} - I)\hat{\phi}||$. Given the definition of $R_{\phi}$ in \eqref{eq:def_phi} and Assumption 4, we can show that 
\begin{subequations}\label{eq:term3_final}
\begin{align}
    ||(R_{\phi}-I)\hat{\phi}|| =\ & || [\hat{\Phi}_y^{\rmt},\hat{\Phi}_u^\rmt]^\rmt (I - \Delta \hat{\Phi}_u)^{-1}\Delta\hat{\phi}_u ||\\
    \leq\ & (||\hat{\Phi}_y|| + ||\hat{\Phi}_u||)\cdot ||(I - \Delta\hat{\Phi}_u)^{-1}||\cdot||\Delta||\cdot||\hat{\phi}_u||\\
    \leq\ & \left(\frac{\rho}{\gamma_1}+\sigma\right)\cdot\frac{\gamma_1}{1-\rho}||\hat{\phi}_u||.
\end{align}
\end{subequations}


Term 4: $||\hat{\Phi}R_{\Phi}(\hat{\Theta} + \tilde{\Theta})||d_W(\overline{\mathbb{P}}_{\mathbf{e}},\mathbb{P}_{\mathbf{e}})$. Considering the definition of $R_{\Phi}$ in \eqref{eq:def_Phi} and applying similar analysis together with Assumption 3 leads to
\begin{equation}\label{eq:term4_final}
    ||\hat{\Phi}R_{\Phi}(\hat{\Theta} + \tilde{\Theta})||d_W(\overline{\mathbb{P}}_{\mathbf{e}},\mathbb{P}_{\mathbf{e}}) \leq \frac{\kappa}{1-\rho}||\hat{\Phi}||\cdot\left( ||\hat{\Theta}||+\gamma_2\right).
\end{equation}

Combining the results in \eqref{eq:term1_final}-\eqref{eq:term4_final} leads to $d_W(\overline{\mathbb{P}}^{\hat{\mathcal{M}}},\mathbb{P}^{\mathcal{M}}) \leq \bar{\varepsilon}$ and constraints \eqref{eq:Wassers_distant}. Finally, for the DRO problem in \eqref{eq:dro_sls} with $\mathbb{Q}\in\mathcal{B}_{\overline{\varepsilon}}(\overline{\mathbb{P}}^{\hat{\mathcal{M}}})$, applying the results of Lemma A.2 in \cite{coulson2021distributionally} by choosing $b=1$, then exploiting Theorem 4.2 in \cite{mohajerin2018data} and considering the SLS parameterization in \eqref{eq:sls} and \eqref{eq:phi_block_lower} with $\mathbf{G}$ replaced by $\mathbf{\hat{G}}$ lead to the optimization problem \eqref{eq:dro_final}. This completes the proof. \hfill $\square$

\remarknew{By fixing the parameters $\rho$ and $\sigma$, which influence the bounds of $||\hat{\Phi}_u||$ and $||\hat{\Phi}_y||$, respectively, the optimization problem \eqref{eq:dro_final} becomes convex and can be efficiently solved. The optimal values of $\rho$ and $\sigma$ can be determined via a grid search by iteratively solving a sequence of convex optimization problems. Alternatively, if $\rho$ and $\sigma$ are treated as decision variables rather than fixed parameters, the problem \eqref{eq:dro_final} becomes nonconvex due to the presence of bilinear terms. Such nonconvexities can still be dealt with using off-the-shelf solvers, such as {\tt Gurobi} and {\tt Ipopt}. Furthermore, the bound derived in \eqref{eq:Wassers_distant} possesses an asymptotic consistency property. In the absence of uncertainties, i.e., when $\gamma_1 = \gamma_2 = \gamma_3 = \kappa = 0$, the last two terms in \eqref{eq:Wassers_distant} vanish. Moreover, the condition $\gamma_1 = 0$ implies that $\rho$ can be chosen arbitrarily small, thereby allowing the Wasserstein distance bound $\bar{\varepsilon}$ to be made arbitrarily close to zero.
}

\remarknew{
In \eqref{eq:dro_sls}, the system constraint is formulated to guarantee that its expected value is within a prescribed bound. To enhance constraint satisfaction and provide probabilistic guarantees, our approach can be readily extended--following the method in \cite{micheli2024data}--to incorporate a conditional value at risk (CVaR) formulation. This extension preserves the structure of the proposed scheme and does not increase computational complexity. Natually, using a CVaR-based formulation to achieve higher probability guarantees results in a more conservative solution.
}

\remarknew{It is worth highlighting that while the parameters $(\mathbf{G},\Theta,\mathbf{y}_0)$ considered in Theorem 1 are related to system matrices $(A,B,C,D)$, the proposed design does not rely on the explicit knowledge of these system matrices. Instead, $(\mathbf{G},\Theta,\mathbf{y}_0)$ can be possibly estimated from system input-output data. As in many existing subspace algorithms (e.g., \cite{van2013closed}), the input-output relationships in \eqref{eq:7} can be leveraged to identify approximate values of $(\mathbf{G},\Theta,\mathbf{y}_0)$. Unlike the case in traditional subspace-based system identification methods, our approach avoids the explicit reconstruction of the system matrices $(A,B,C,D)$, a task that is typically complex and challenging.
}

\section{Simulation Results}\label{sec:simu}
In this section, numerical simulation results are presented to demonstrate the effectiveness of the proposed approach. The numerical study is based on the linear system considered in \cite{breschi2023data}, with a more challenging modification by setting $B_2 = 0$ such that small perturbations of $B_2$ will change the sign of the control input, resulting in different behaviours. The true system dynamics used in simulation are 
\begin{subequations}\label{eq:simu_true}
\begin{align}
    &x_{t+1} = \begin{bmatrix}
        0.7326 & -0.0861\\
        0.1722 & 0.9909
    \end{bmatrix}x_{t} +
    \begin{bmatrix}
        0.0609\\
        0
    \end{bmatrix}u_t + w_t,\\
    &y_t = \begin{bmatrix}
        0 & 1.4142
    \end{bmatrix}x_t + v_t,
\end{align}
\end{subequations}
where $w_t\in\mathbb{R}^2$ and $v_t\in\mathbb{R}$ are assumed to follow uniform distributions with each element belonging to $\mathcal{U}(-0.01,0.01)$. For constructing the innovation process, the feedback matrix $L$ in \eqref{eq:sys_dyn} is set as $[0.1, 0.1]^{\mathrm{T}}$, and $A - LC$ is Schur stable. For the SLS design, its prediction horizon is set as $T = 15$, and $\tau = 25$. The cost function is set as
\begin{equation}
    h(\bm{\eta}) := \left\lVert
    \begin{bmatrix}
    \mathbf{Q} & \\
     & \mathbf{R}
    \end{bmatrix}
    \begin{bmatrix}
        \mathbf{y}\\
        \mathbf{u}
    \end{bmatrix}
    \right\rVert,
\end{equation}
where $\mathbf{Q}$ and $\mathbf{R}$ are diagonal matrices with diagonal elements as $1$ and $0.1$, respectively. System constraints at each time step $k=0,\cdots,T$ are defined as 
\begin{equation*}
    g_k(\bm{\eta}):=
    \max\{-y_k - 0.01, u_k - 1, -u_k - 1\} <= 0.
\end{equation*}
This constraint limits the output and input signals satisfying $y_t \geq -0.01$ and $-1\leq u_t\leq 1$, respectively.

The upper bounds of uncertainties are $\gamma_1 = \gamma_2= \gamma_3 = 0.01$. In the simulation, 50 samples of system matrices $(\hat{A},\hat{B},\hat{C})$ are randomly generated such that the corresponding parameters $(\hat{\mathbf{G}}, \hat{\Theta}, \hat{\mathbf{y}}_0)$ reside in the prescribed scope. Random input signals within $[-1,1]$ are used to generate 100 samples of the innovation process signal $\mathbf{e}$. The upper bound of the Wasserstein distance $d_W(\overline{\mathbb{P}}_{\mathbf{e}},\mathbb{P}_{\mathbf{e}})$ is set as $\kappa = 0.005$. 

In our simulation, we mainly investigate two SLS approaches:
\begin{itemize}
    \item N-SLS: nominal system level synthesis design adopting the CE and SAA solutions, where the nominal parameters $(\hat{\mathbf{G}}, \hat{\Theta}, \hat{\mathbf{y}}_0)$ and the averaged historical samples of $\mathbf{e}$ are used in Theorem 1 of Section \ref{sec:sls_output_affine}. 
    \item DR-SLS: distributionally robust system level synthesis design proposed in this work by solving \eqref{eq:dro_final}.
\end{itemize}
For each sample of the parameters $(\hat{\mathbf{G}}, \hat{\Theta}, \hat{\mathbf{y}}_0)$, after solving the SLS design problems with the above two approaches, the control policy parameters $(\mathbf{K},\mathbf{p})$ are used for performing the closed-loop simulation with the true system dynamics in \eqref{eq:simu_true}. 

\begin{figure}[htb]
    \centering
    \includegraphics[width=\linewidth]{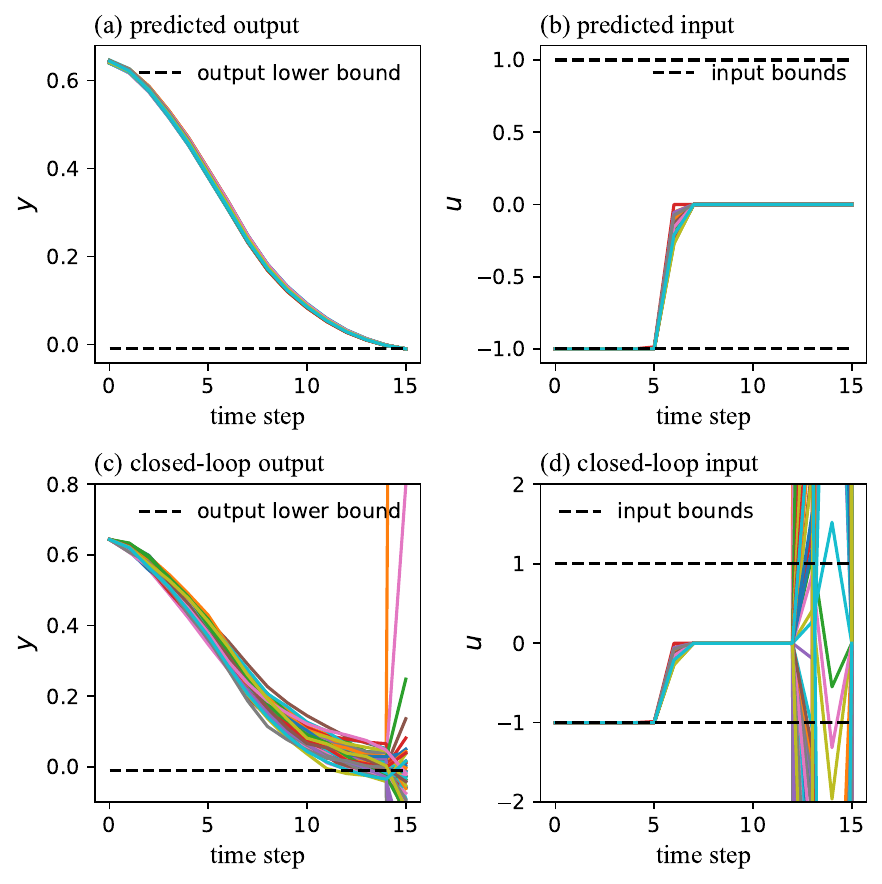}
    \caption{Output and input trajectories with N-SLS approach: (a) open-loop output prediction, (b) open-loop input prediction, (c) closed-loop output, (d) closed-loop input.}
    \label{fig:nsls_simu}
\end{figure}

\begin{figure}[htb]
    \centering
    \includegraphics[width=\linewidth]{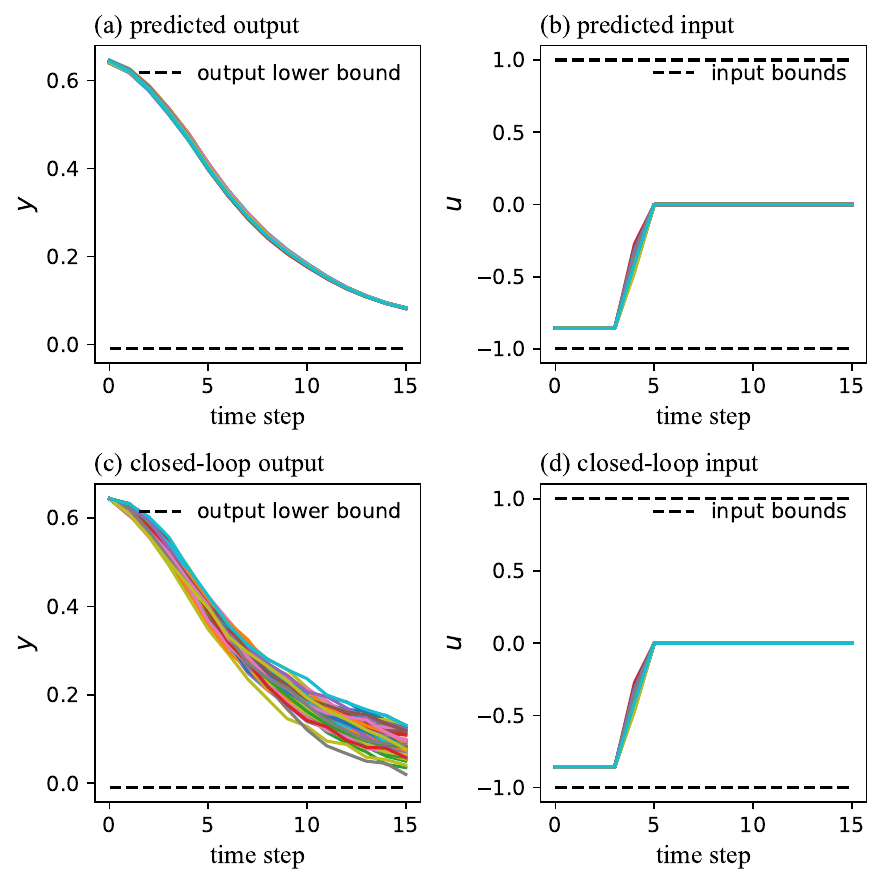}
    \caption{Output and input trajectories with DR-SLS approach: (a) open-loop output prediction, (b) open-loop input prediction, (c) closed-loop output, (d) closed-loop input.}
    \label{fig:drsls_simu}
\end{figure}

\begin{figure}[htb]
    \centering
    \includegraphics[width=\linewidth]{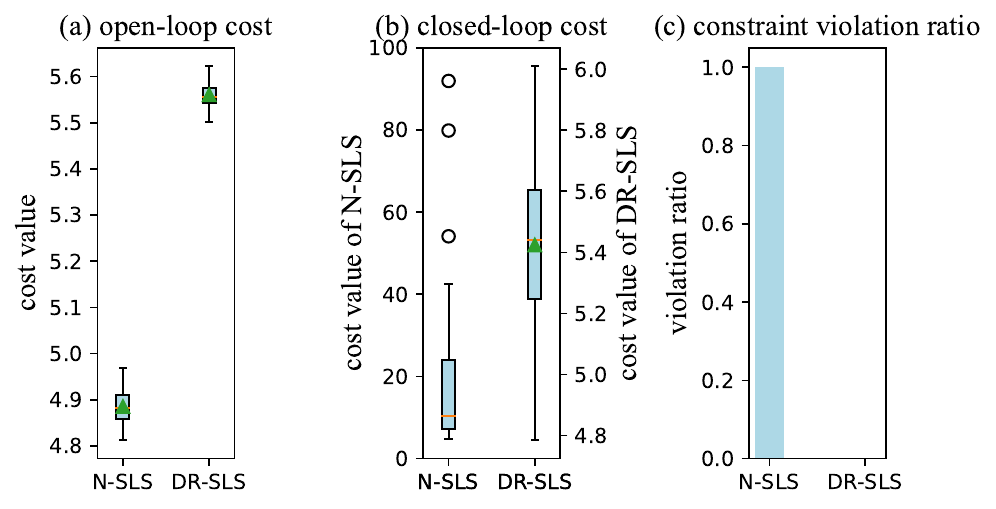}
    \caption{Cost values and ratios of constraint violation of N-SLS and DR-SLS: (a) open-loop cost values, (b) closed-loop cost values, (c) ratios of constraint violation. }
    \label{fig:obj_simu}
\end{figure}

Simulation results are shown in Figs.~\ref{fig:nsls_simu}--\ref{fig:obj_simu}. Fig.~\ref{fig:nsls_simu} depicts the output and input trajectories adopting the N-SLS approach for both open-loop prediction and closed-loop simulation. It can be observed that the N-SLS approach is aggressive in optimizing the cost function such that the output trajectory is close to the predefined lower bound. However, while the predicted output and input profiles with the approximated system parameters satisfy all system constraints and converge to the expected setpoints, shown in Fig.~\ref{fig:nsls_simu}(a-b), the closed-loop response with the true model becomes unstable and the output-input trajectories diverge and violate system constraints, see Fig.~\ref{fig:nsls_simu}(c-d).

Fig.~\ref{fig:drsls_simu} presents the simulation results adopting the proposed DR-SLS approach. It can be observed from Fig.~\ref{fig:drsls_simu}(a-b) that the output and input predictions all satisfy system constraints and leave margins to the respective permissible boundaries. In comparison with the case of the N-SLS approach shown in Fig.~\ref{fig:nsls_simu}, the input-output predictions do not exactly converge to the origin, and small offsets are introduced to enhance robustness against possible model uncertainties. Fig.~\ref{fig:drsls_simu}(c-d) show the true closed-loop trajectories. In contrast to the unstable behavior observed with the N-SLS approach, the DR-SLS approach ensures closed-loop stability while satisfying all system constraints, demonstrating its robustness and effectiveness.

Fig.~\ref{fig:obj_simu} compares the cost values and constraint violation ratios of the two approaches for all samples of model parameters. It is clear from Fig.~\ref{fig:obj_simu}(a) that the predictive cost values of the N-SLS approach are lower than those of the DR-SLS approach. However, as shown in Fig.~\ref{fig:obj_simu}(b), when implementing the corresponding control policies in closed-loop simulation, the closed-loop cost value of the N-SLS approach is tremendously larger than the DR-SLS approach since the closed-loop responses with the control policy computed via the N-SLS approach become divergent. In addition, Fig.~\ref{fig:obj_simu}(c) shows that the N-SLS approach results in constraint violations in all closed-loop simulations, whereas the proposed DR-SLS approach robustly enforces constraint satisfaction in the presence of model mismatch and additive disturbances.

In summary, compared to the nominal SLS approach that adopts SAA and CE solutions to handle model uncertainties, the proposed DR-SLS formulation provides a systematic and principled method for enhancing robustness. Naturally, this improved robustness comes at the cost of increased conservatism in the resulting optimal solution.

\section{Conclusion}\label{sec:conclusion}
In this paper, we presented a distributionally robust system level synthesis design with output feedback affine control policy for uncertain linear systems subject to additive model mismatch and stochastic disturbances. By designing the SLS parameterization for the linear system in its innovation form, the proposed SLS scheme, compared to existing methods, reduces the number of decision variables and constraints and accommodates a broader class of system structures. The DR-SLS problem is cast as a robust optimization task that optimizes the expectation of the cost function while ensuring constraint satisfaction under the worst-case uncertainty distribution. We analyzed the impact of model mismatch and the empirical estimation error of stochastic uncertainties on the distribution shift between predictive closed-loop response and real closed-loop response, and derived an upper bound in the sense of the Wasserstein metric. Leveraging tools from robust SLS and DRO, we developed a tractable convex reformulation for the proposed DR-SLS problem. Numerical results demonstrate that the proposed scheme ensures safe and effective control despite the model mismatch and uncertainties. 

\section*{References}
\bibliographystyle{IEEEtran.bst}
\bibliography{ref}

\end{document}